\documentclass[12pt]{article}
\usepackage[a4paper,margin=1in]{geometry}
\usepackage[T1]{fontenc}
\usepackage{lmodern}
\usepackage{microtype}
\usepackage{amsmath,amssymb,amsthm,mathtools}
\usepackage[hidelinks]{hyperref}
\usepackage[
    backend=biber,
    style=numeric,
    sorting=nyt,
    giveninits=true,
    doi=true,
    url=false,
    isbn=false
]{biblatex}
\newtheorem{maintheorem}{Theorem}
\newtheorem{theorem}{Theorem}[section]
\newtheorem{proposition}[theorem]{Proposition}
\newtheorem{lemma}[theorem]{Lemma}
\newtheorem{corollary}[theorem]{Corollary}
\theoremstyle{definition}
\newtheorem{definition}[theorem]{Definition}

\theoremstyle{plain}
\newtheorem{theoremA}{Theorem}

\numberwithin{equation}{section}
\newcommand{\D}{\mathbb D}
\newcommand{\N}{\mathbb N}

\newcommand{\dist}{\operatorname{dist}}
\newcommand{\1}{\mathbf 1}
\newcommand{\firstpagefunding}[1]{%
    \begingroup
    \renewcommand{\thefootnote}{}%
    \begin{NoHyper}%
    \footnotetext{#1}%
    \end{NoHyper}%
    \addtocounter{footnote}{-1}%
    \endgroup
}

\begin{document}
\begin{center}
	{\scshape\Large Finite Angular Derivative Sets
		and Almost Isometries of the Hyperbolic Disc\par}
	\vspace{0.45cm}
	{\Large Angel Lorenzo \par}
	\vspace{0.35cm}
\end{center}

\firstpagefunding{The author is supported in part by the Spanish Ministerio de Ciencia, Innovación y Universidades (project PID2025-167726NB-I00).}
\begin{abstract}
We prove that a subset $E$ of the unit circle is the set of points at which an analytic self-map of the unit disc has a finite angular derivative if and only if $E$ is an $F_\sigma$ set. We also prove that, for every countable dense subset $E$ of the unit circle, there exists an almost isometry in the sense of Garnett--Papadimitrakis whose angular derivative is finite precisely on $E$.
\end{abstract}
\section{Introduction}

Let $\D=\{z\in\mathbb C:|z|<1\}$ be the unit disc and let $\partial\D$ be the unit circle. Let $\mathcal S$ denote the \emph{Schur class}, that is, the class of analytic mappings $f:\D\to\D$. A function $f\in\mathcal S$ has a \emph{finite angular derivative} at a point $\zeta\in\partial\D$ if the non-tangential limit of $f(z)$ is a point $\eta\in\partial\D$, and $\frac{f(z)-\eta}{z-\zeta}$ has a finite non-tangential limit at $\zeta$, which is denoted by $f'(\zeta)$. By the Julia--Carathéodory theorem (see, for instance,
\cite[Chapter~4]{Shapiro}), this is equivalent to
\begin{equation}\label{eq:JC-intro}
	\liminf_{z\to\zeta}
	\frac{1-|f(z)|}{1-|z|}<\infty.
\end{equation}
When \eqref{eq:JC-intro} fails, it is common to denote
$|f'(\zeta)|=\infty$. Angular derivatives were introduced by Carathéodory and have become a classical notion in complex analysis; see for instance
\cite{CMR}, \cite{Mashreghi}, and \cite{Shapiro}.
For $f\in\mathcal S$, write
$\mathcal F(f)=\{\zeta\in\partial\D: f\text{ has a finite angular derivative at }\zeta\}$.
Our first result gives a topological description of the sets that
can occur in this way.

\begin{maintheorem}\label{thm:intro-fsigma}
	\leavevmode\par
	\begin{enumerate}
		\item[(a)] If $f\in\mathcal S$, then $\mathcal F(f)$ is an
		      $F_\sigma$ subset of $\partial\D$.
		\item[(b)] If $E\subset\partial\D$ is an $F_\sigma$ set, then there exists a Blaschke product $B$ such that $\mathcal F(B)=E$.
	\end{enumerate}
\end{maintheorem}

Part~(a) follows easily from well-known results on angular derivatives. Part~(b) is harder and uses a dyadic construction: one places the zeros of a Blaschke product so that its angular derivative is finite precisely on the prescribed set $E$.

A more rigid problem arises when one considers subclasses of the Schur class $\mathcal{S}$. Let
$d_h(z,w)$ be the hyperbolic distance between the points $z,w\in\D$. For a set $A\subset\D$, let $\operatorname{diam}_h(A)$ denote
its hyperbolic diameter. Given $c>0$,
let $M(c)$ be the class of functions $f\in\mathcal S$ satisfying
\begin{equation}\label{eq:Mc-intro}
	\operatorname{diam}_h f(\mathcal B)
	\geq \operatorname{diam}_h(\mathcal B)-c
\end{equation}
for every hyperbolic disc $\mathcal B\subset\D$. A function $f\in\mathcal{S}$ is called an almost isometry if $f\in M(c)$ for some $c>0$. Garnett and Papadimitrakis \cite{GP91} proved that an almost isometry must be a Blaschke product, and gave characterizations of the Blaschke products that are almost isometries in terms of angular derivatives and the distribution of their zeros. Geometric classes of analytic self-maps, defined through almost preservation of hyperbolic area and bounded compression, have recently been studied in \cite{IvriiNicolau2024} and \cite{IvriiNicolau2025} using angular derivatives. One consequence of the Garnett--Papadimitrakis results is that $\mathcal F(B)$ must be dense in $\partial\D$ whenever $B$ is an almost isometry. In the same
paper, they include an example, due to Carleson and Jones, of
an almost isometry for which the set of points admitting a finite angular
derivative has Lebesgue measure zero. This raises the natural question of how small $\mathcal F(B)$ can be for an almost isometry $B$. Our second result shows that no size requirement other than density is needed.
\begin{maintheorem}\label{thm:intro-countable}
	Let $E\subset\partial\D$ be a countable dense set. Then there exists an almost isometry $B$ such that $\mathcal F(B)=E$.
\end{maintheorem}

The proof of Theorem~\ref{thm:intro-countable} is constructive and is based on a dyadic analysis designed to satisfy simultaneously two competing conditions: the zeros of $B$ must be sufficiently well distributed to guarantee the almost isometry estimates, while their approach to the unit circle must be arranged so that the angular derivative is finite exactly at the prescribed set $E$.

The paper is organized as follows. Section~2 collects preliminary results on
angular derivatives, Blaschke products, and almost isometries. In Section~3, Theorem~\ref{thm:intro-fsigma} is proved. We also present a result analogous to Theorem~\ref{thm:intro-fsigma} for outer functions in $\mathcal{S}$, where an additional measure-theoretic obstruction appears.
Section~4 introduces the dyadic construction underlying the proof of
Theorem~\ref{thm:intro-countable}. Section~5 is devoted to the proof of Theorem~\ref{thm:intro-countable}, where we first establish the result when $E$ consists of the dyadic points of $\partial\D$ and then extend it to arbitrary countable dense subsets of the unit circle.

This work is part of the author's PhD thesis at Universitat Autònoma de Barcelona.
\paragraph{AI Declaration.}
ChatGPT (OpenAI, GPT-5.6 Luna) was used to correct typographical and grammatical errors and to improve the conciseness of the introduction.

\section{Preliminaries}
\label{sec:preliminaries}
Let $\{a_n\}$ be a sequence of points of $\D$ satisfying $\sum_n(1-|a_n|)<+\infty$. The Blaschke product with zeros $\{a_n\}$ is, with the usual convention when $a_n=0$,
\begin{equation}\label{eq:blaschke-product}
	B(z)=\prod_n \frac{|a_n|}{a_n}\,
	\frac{a_n-z}{1-\overline{a_n}z}, \ \ z\in\D.
\end{equation}
Given a positive Borel measure, singular with respect to Lebesgue measure $dm$ on $\partial\D$, let $S_\mu$ be the singular inner function defined as
\begin{equation}\label{eq:singular-inner}
	S_\mu(z)=
	\exp\!\left(-\int_{\partial\D}
	\frac{\xi+z}{\xi-z}\,d\mu(\xi)\right), \ \ z\in\D.
\end{equation}
A function $O\in\mathcal{S}$ is \emph{outer} if it can be written
as
\begin{equation}\label{eq:outer-factor}
	O(z)=e^{i\gamma}
	\exp\!\left(-\int_{\partial\D}
	\frac{\xi+z}{\xi-z}\,q(\xi)\,dm(\xi)\right), \ \ z\in\D,
\end{equation}
where $q\in L^1(\partial\D)$, $q\ge0$, and $e^{i\gamma}\in\partial\D$. For later use, define the boundary potential
\begin{equation}\label{eq:outer-potential}
	U_q(\zeta)=\int_{\partial\D}\frac{q(\xi)}{|\xi-\zeta|^2}\,dm(\xi), \qquad \zeta\in\partial\D.
\end{equation}
The classical inner-outer factorization says that for any nonzero $f\in \mathcal{S}$, there exists a sequence of points $\{a_n\}\subset\D$, $\sum_n(1-|a_n|)<+\infty$, a positive singular measure $\mu$, a nonnegative function $q\in L^1(\partial\D)$ and $e^{i\gamma}\in\partial\D$ such that $f=BS_\mu O$. See, for instance, \cite{Duren}.
The following standard form of the angular-derivative criterion combines these three factors (see Cima--Matheson--Ross
\cite[Section~1.7]{CMR} or \cite{AhernClark}).
\begin{lemma}[\cite{CMR}]\label{lem:canonical-factor}
	Let $f=BS_\mu O\in\mathcal S$ be nonzero, with factors given by (\ref{eq:blaschke-product}), (\ref{eq:singular-inner}) and (\ref{eq:outer-factor}). Then $f$ has a finite angular derivative at
	$\zeta\in\partial\D$ if and only if
	\begin{equation}\label{eq:canonical-angular-criterion}
		\sum_n\frac{1-|a_n|^2}{|\zeta-a_n|^2}
		+2\int_{\partial\D}\frac{d\mu(\xi)}{|\xi-\zeta|^2}
		+2\int_{\partial\D}
		\frac{q(\xi)}{|\xi-\zeta|^2}\,dm(\xi)
		<\infty.
	\end{equation}
	Whenever this condition holds, one has
	\begin{equation}\label{eq:canonical-angular-formula}
		|f'(\zeta)|=
		\sum_n\frac{1-|a_n|^2}{|\zeta-a_n|^2}
		+2\int_{\partial\D}\frac{d\mu(\xi)}{|\xi-\zeta|^2}
		+2\int_{\partial\D}
		\frac{q(\xi)}{|\xi-\zeta|^2}\,dm(\xi).
	\end{equation}
	In particular,
	\begin{align}
		\mathcal F(B)
		 & =\left\{\zeta\in\partial\D:
		\sum_n\frac{1-|a_n|^2}{|\zeta-a_n|^2}<\infty\right\},
		\label{eq:calF-Blaschke}                        \\
		\mathcal F(S_\mu)
		 & =\left\{\zeta\in\partial\D:
		\int_{\partial\D}\frac{d\mu(\xi)}{|\xi-\zeta|^2}<\infty\right\},
		\label{eq:calF-singular}                        \\
		\mathcal F(O)
		 & =\left\{\zeta\in\partial\D:\int_{\partial\D}
		\frac{q(\xi)}{|\xi-\zeta|^2}\,dm(\xi)<\infty\right\},
		\label{eq:calF-outer}
	\end{align}
	and hence
	\begin{equation}\label{eq:calF-intersection}
		\mathcal F(f)=\mathcal F(B)\cap \mathcal F(S_\mu)\cap \mathcal F(O).
	\end{equation}
\end{lemma}

Thus the three canonical factors contribute independently to the angular derivative.
\begin{lemma}\label{lem:potential-lsc}
	Let $f\in\mathcal{S}$ be as in Lemma~\ref{lem:canonical-factor}. Define the function
	$|f'|:\partial\D\to[0,+\infty]$ by the right-hand side of
	equation (\ref{eq:canonical-angular-formula}). Then $|f'|$ is lower
	semicontinuous on $\partial\D$.
\end{lemma}

\begin{proof}
	Write $f=BS_\mu O$. We show that each of the three terms in
	equation (\ref{eq:canonical-angular-formula}) is lower semicontinuous. The first term
	$\sum_{n}\frac{1-|a_n|^2}{|\zeta-a_n|^2}$ is lower semicontinuous on $\partial\D$ because it is an increasing pointwise limit of continuous functions on $\partial\D$. For $M>0$, put $K_M(\xi,\zeta)=\min\{M,|\xi-\zeta|^{-2}\}$, where $K_M(\xi,\xi)=M$, and define
	\begin{equation}
		V_{\mu,M}(\zeta)
		=
		\int_{\partial\D}K_M(\xi,\zeta)\,d\mu(\xi),
		\qquad
		U_{q,M}(\zeta)
		=
		\int_{\partial\D}K_M(\xi,\zeta)q(\xi)\,dm(\xi).
	\end{equation}
	Since the kernel $K_M$ is continuous and bounded, dominated convergence
	shows that both $V_{\mu,M}$ and $U_{q,M}$ are continuous on
	$\partial\D$. Monotone convergence gives
	\begin{equation}
		\int_{\partial\D}\frac{d\mu(\xi)}{|\xi-\zeta|^2}
		=
		\sup_{M>0}V_{\mu,M}(\zeta),
		\qquad
		U_q(\zeta)
		=
		\sup_{M>0}U_{q,M}(\zeta).
	\end{equation}
	Therefore both integral terms are lower semicontinuous. Thus, $|f'|$ is lower semicontinuous on $\partial\D$.
\end{proof}
We now present a result of Garnett and Papadimitrakis which we will use in the proof of Theorem~\ref{thm:intro-countable}. For a Blaschke sequence $Z\subset\D$, define
\begin{equation}\label{eq:SZ}
	S_Z(\xi)=\sum_{a\in Z}\frac{1-|a|^2}{|\xi-a|^2}, \qquad \xi\in\partial\D.
\end{equation}
If $I\subset\partial\D$ is an arc with normalized length $m(I)$ and centre $c_I\in\partial\D$, set $z_I=(1-m(I))c_I$ and
\begin{equation}\label{eq:TZ}
	T_Z(I)=\sum_{a\in Z}\frac{1-|a|^2}{|1-\overline a z_I|^2}.
\end{equation}

\begin{theoremA}[Garnett--Papadimitrakis]\label{thm:GP-criterion}
	Let $B$ be a Blaschke product with zero set $Z$. Then $B$ is an almost isometry if and only if there is a constant $A>0$ such that, for every arc $I\subset\partial\D$ with $0<m(I)<1/2$, the following two conditions hold:
	\begin{equation*}\tag{AI1}\label{AI1}
		\inf_{\xi\in I} S_Z(\xi)\le A\,T_Z(I),
	\end{equation*}
	and
	\begin{equation*}\tag{AI2}\label{AI2}
		(1-|z_I|^2)T_Z(I)\le A.
	\end{equation*}
\end{theoremA}

Condition \eqref{AI2} means that the measure $\nu_Z=\sum_{a\in Z}(1-|a|^2)\delta_a$ is a Carleson measure, that is, there exists a constant $C>0$ such that $\nu_Z(Q(I))\leq C m(I)$ for every arc $I\subset\partial\D$.  Here, $Q(I):=\{r\xi: \xi\in I, \ 1-m(I)<r<1\}$ is a Carleson box based on the arc $I\subset\partial\D$. See, for example, \cite[Chapter~VI, \S3]{Garnett}. Note that \eqref{AI1} implies that $\mathcal F(B)$ is dense in $\partial\D$ if $B$ is an almost isometry.
\section{Finite angular derivative sets}
\label{sec:finite-sets}
For part (b) of Theorem~\ref{thm:intro-fsigma}, we will use
the following two elementary results. Roughly speaking, given a compact set $K\subset\partial\D$ we will construct a finite set of zeros while keeping both its Blaschke mass and its influence away from a prescribed open neighbourhood of $K$ arbitrarily small.

\begin{lemma}[Local block]\label{lem:local-block}
	Let $K\subset U\subset\partial\D$, with $K$ compact and $U$ open. Given
	$\varepsilon>0$ and $c>0$, there is a finite multiset $Z\subset\D$
	such that
	\begin{equation}\label{eq:local-block}
		\sum_{a\in Z}(1-|a|)\le\varepsilon,
		\qquad \sup_{\xi\in \partial\D\setminus U}S_Z(\xi)\le\varepsilon,
		\qquad \inf_{\xi\in K} S_Z(\xi)\ge c.
	\end{equation}
\end{lemma}

\begin{proof}
	Choose finitely many pairs of arcs $(J_i,I_i)$, $i=1,\ldots,M$, such that $K\subset \bigcup_{i=1}^M J_i\subseteq \bigcup_{i=1}^M \overline{J_i}\subseteq\bigcup_{i=1}^M I_i\subset U$. On one
	such $J=\{e^{it}:\alpha\le t\le\beta\}$, put
	$\delta=\dist(J,\partial\D\setminus I)>0$. Choose $N$ equally spaced angular points $e^{i\theta_j}$ on $J$, that is, with separation
	$h\asymp m(J)/N$, and set $\rho=4ch^2$, $a_j=(1-\rho)e^{i\theta_j}$, and $Z_{J,I}=\{a_1,\ldots,a_N\}$.
	For large $N$, $\rho\le h$ and $\rho\le\delta/2$. If
	$\xi=e^{it}\in J$, choose $j$ with $|t-\theta_j|\le h$. Then
	$|\xi-a_j|\le h+\rho\le2h$ and $1-|a_j|^2\ge\rho$, so this single zero contributes at least $\rho/(4h^2)=c$.
	Meanwhile, $\sum_{j=1}^N(1-|a_j|)=N\rho\asymp cm(J)^2/N\longrightarrow0$.
	For $\xi\notin I$, we have $|\xi-a_j|\gtrsim\delta$, uniformly in
	$j$, and hence $S_{Z_{J,I}}(\xi)\lesssim N\rho/\delta^2\longrightarrow0$.
	Thus both upper estimates may be made arbitrarily small without changing the
	lower bound on $J$. For the general compact set, choose finitely many pairs $(J_r,I_r)$ such that $\overline{J_r}\subset I_r\subset \overline{I_r}\subset U$ whose interiors cover $K$, apply the arc
	construction with error $\varepsilon/M$, and take the multiset union.
\end{proof}
The following auxiliary result is standard. Its short proof is included for the sake of completeness.
\begin{lemma}\label{lem:layers}
	For every non-empty proper open set $U\subset\partial\D$, there exist compact sets $K_r$ and open sets
	$V_r$ of $\partial\D$ such that $K_r\subset V_r\subset U$, $U=\bigcup_{r\ge0}K_r$, and
	$\sum_{r\ge0}\1_{V_r}\le4$.
\end{lemma}

\begin{proof}
	Let $d(\xi)=\dist(\xi,\partial\D\setminus U)$. Take $K_0=\{\xi\in \partial\D:d(\xi)\ge1/2\}$, $V_0=\{\xi\in\partial\D:d(\xi)>1/4\}$,
	and, for $r\ge1$, $K_r=\{\xi\in\partial\D:2^{-(r+1)}\le d(\xi)\le2^{-(r-1)}\}$, $V_r=\{\xi\in\partial\D:2^{-(r+2)}<d(\xi)<2^{-(r-2)}\}$.
	Continuity of $d$ gives compactness, openness, and the stated inclusions.
	The $K_r$ cover $U=\{\xi\in\partial\D:d(\xi)>0\}$. For a fixed $\xi\in U$, membership in
	$V_r$ restricts $r$ to an interval of length four on the logarithmic
	scale, so at most four indices occur. Empty layers may be discarded.
\end{proof}

\begin{proof}[Proof of Theorem~\ref{thm:intro-fsigma}] \textbf{Necessity.} By Lemma~\ref{lem:potential-lsc}, the function
    $\zeta\mapsto |f'(\zeta)|$ is lower semicontinuous on $\partial\D$. Hence, for every $N\geq 1$, the set
    $\{\zeta\in\partial\D:|f'(\zeta)|\leq N\}$ is closed. Moreover,
\begin{equation}
	\mathcal F(f)
	=
	\bigcup_{N\geq 1}
	\{\zeta\in\partial\D:|f'(\zeta)|\leq N\}.
\end{equation}
Therefore $\mathcal F(f)$ is an $F_\sigma$ set.

	\medskip\noindent\textbf{Sufficiency.}
	Let $F\subset\partial\D$ be an $F_\sigma$ set, and put $E=\partial\D\setminus F$. If $E=\varnothing$, take a finite Blaschke product. If $E=\partial\D$, apply Lemma~\ref{lem:local-block} globally at level $n$ with mass $2^{-n}$ and lower bound one; by \eqref{eq:calF-Blaschke}, the resulting Blaschke product satisfies $\mathcal F(B)=\varnothing=F$. We may therefore assume that $E$ is a nonempty proper $G_\delta$ set. Write $E=\bigcap_{n\ge1}U_n$, where $\{U_n\}$ is a decreasing sequence of proper open subsets of $\partial\D$. For each $n$, apply Lemma~\ref{lem:layers} to obtain
	$K_{n,r}\subset V_{n,r}\subset U_n$, $r=1,2,\ldots$. Apply Lemma~\ref{lem:local-block} to each pair with $c=1$ and errors
	$\varepsilon=\varepsilon_{n,r}>0$ satisfying $\sum_{n,r}\varepsilon_{n,r}<\infty$.
	Write $Z_{n,r}$ for the resulting finite collection of points, $S_{n,r}$ for its boundary
	sum, and let $Z=\bigcup_{n,r}Z_{n,r}$. Then
	$\sum_{a\in Z}(1-|a|)\le\sum_{n,r}\varepsilon_{n,r}<\infty$, so $Z$ is a Blaschke sequence.

	If $\xi\in E$, then for any $n$ there exists $r(n)$ such that $\xi\in K_{n,r(n)}$. Therefore
	\begin{equation}\label{eq:suff-divergence}
		\sum_{a\in Z}\frac{1-|a|^2}{|\xi-a|^2}
		\ge\sum_{n\ge1}S_{n,r(n)}(\xi)
		\ge\sum_{n\ge1}1=\infty.
	\end{equation}

	If $\xi\notin E$, choose $N$ such that $\xi\notin U_n$ for every
	$n\ge N$. Since $V_{n,r}\subset U_n$, the tail satisfies
	$\sum_{n\ge N}\sum_rS_{n,r}(\xi)\le\sum_{n\ge N}\sum_r\varepsilon_{n,r}<\infty$.
	For a fixed $n<N$, the bounded-overlap property of Lemma~\ref{lem:layers} leaves only finitely many
	indices $r$ for which $\xi\in V_{n,r}$. Each of the corresponding
	blocks is finite. For every other $r$, the exterior estimate gives
	$S_{n,r}(\xi)\le\varepsilon_{n,r}$. Hence the initial contribution is finite as well, and the sum diverges exactly on $E$.
\end{proof}
\subsection{Outer functions}
The analogue of Theorem~\ref{thm:intro-fsigma} for outer functions requires an additional condition.

\begin{theorem}\label{thm:outer}
	Let $E\subsetneq\partial\D$ and $D=\partial\D\setminus E$. There exists an outer function $O\in\mathcal S$ with $\mathcal F(O)=E$ if and only if $E$ is $F_\sigma$ and $m(I\cap D)>0$ for every open arc $I\subset\partial\D$ such that $I\cap D\neq\varnothing$.
\end{theorem}

\begin{proof}
	Suppose first that $\mathcal F(O)=E$, and let
	$q=-\log|O^*|$. By Lemma~\ref{lem:canonical-factor},
	$E=\{\zeta\in\partial\D:U_q(\zeta)<\infty\}$, which is an $F_\sigma$ set. Let $I_r(\zeta)$ be the arc of $\partial\D$ centred at $\zeta$ with $m(I_r(\zeta))=r$. If $U_q(\zeta)<\infty$, then
	\begin{equation}\label{eq:outer-averages}
		\frac{1}{m(I_r(\zeta))}
		\int_{I_r(\zeta)}q\,dm
		\lesssim rU_q(\zeta)
		\longrightarrow0, \ \ \text{as} \ r\to 0.
	\end{equation}
	At every Lebesgue point $\zeta\in E$ of $q$, the left-hand side converges to $q(\zeta)$. Hence $q(\xi)=0$ for almost every $\xi\in E$.

	Suppose now that an open arc $I$ satisfies $I\cap D\neq\varnothing$ but $m(I\cap D)=0$. Then $q=0$ almost everywhere on $I$. Choose $x\in I\cap D$ and an arc $J$ with
	$x\in J\subset\overline J\subset I$. Therefore
	$U_q(x)\le \dist(x,\partial\D\setminus J)^{-2}\|q\|_1<\infty$, contradicting $x\in D$. This proves the necessity of $m(I\cap D)>0$.

	Conversely, if $E=\varnothing$, take $q\equiv1$. Then $U_q(\zeta)=\infty$ for every $\zeta\in\partial\D$, and the corresponding outer function is the constant $e^{-1}$. Assume henceforth that $E\neq\varnothing$. Write $E=\bigcup_{n\ge1}K_n$, where $\{K_n\}$ is an increasing sequence of compact subsets of $\partial\D$, and let $G_n=\partial\D\setminus K_n$. For each $n$, take a Whitney decomposition
	$\mathcal W_n=\{J_{k,n}\}_k$ of $G_n$. Thus
	$G_n=\bigcup_k J_{k,n}$, the arcs $J_{k,n}$ have pairwise disjoint interiors, and there exists a universal constant $c>0$ such that
	\begin{equation}\label{eq:whitney}
		c^{-1}  m(J_{k,n})\leq \dist(J_{k,n},K_n)\leq c  m(J_{k,n}), \ \ k,n=1,2,\ldots.
	\end{equation}
	Choose $\varepsilon_n>0$ with $\sum_n\varepsilon_n<\infty$. For each Whitney arc $J_{k,n}$, choose a finite family $B(k,n)$ of arcs $J\subset J_{k,n}$, with pairwise disjoint interiors, covering $J_{k,n}$, and satisfying
	$m(J)\asymp\varepsilon_n m(J_{k,n})^2$ for any $J\in B(k,n)$. Then
	\begin{equation}\label{eq:Jsum}
		\sum_{J\in B(k,n)}m(J)^2
		\lesssim
		\varepsilon_n m(J_{k,n})^3.
	\end{equation}

	If $J\in B(k,n)$ satisfies $J\cap D\neq\varnothing$, by assumption we have $m(J\cap D)>0$. Choose a compact set
	$A_{n,J}\subset J\cap D$ of positive measure, and define
	\begin{equation}\label{eq:outer-density-pieces}
		q_{n,J}
		:=
		\frac{m(J)^2}{m(A_{n,J})}\1_{A_{n,J}},
		\qquad
		q_n
		:=
		\sum_k\sum_{J\in B(k,n)}q_{n,J},
		\qquad
		q
		:=
		\sum_{n\ge1}q_n.
	\end{equation}
	For $J\in B(k,n)$ with $J\cap D=\varnothing$, set $q_{n,J}=0$. By \eqref{eq:Jsum},
	\begin{equation}\label{eq:outer-qn-mass}
		\|q_n\|_1
		\le
		\sum_k\sum_{J\in B(k,n)}m(J)^2
		\lesssim
		\varepsilon_n\sum_k m(J_{k,n})^3
		\lesssim
		\varepsilon_n\sum_k m(J_{k,n})
		\lesssim
		\varepsilon_n.
	\end{equation}
	Hence Tonelli's theorem gives $q\in L^1(\partial\D)$, and by construction $q\ge0$ almost everywhere. At the same time, if $x\in J$, then every
	$\xi\in A_{n,J}\subset J$ satisfies
	$|\xi-x|\lesssim m(J)$, and therefore
	\begin{equation}\label{eq:outer-piece-lower}
		U_{q_{n,J}}(x)=\int_{\partial\D}\frac{q_{n,J}(\xi)}{|x-\xi|^2}dm(\xi)
		\gtrsim1.
	\end{equation}
	If $x\in D$, then $x\in G_n$ for every $n$. Hence, for every $n$, there are indices $k$ and $J\in B(k,n)$ such that $x\in J$. Since $x\in D$, \eqref{eq:outer-piece-lower} gives $U_q(x)=\infty$.

	Assume now that $\zeta\in K_N$. Since the sequence $\{K_n\}$ is increasing, $\zeta\in K_n$ for all $n\ge N$. If $J\in B(k,n)$, then $A_{n,J}\subset J\subset J_{k,n}$, and \eqref{eq:whitney} gives
	$\dist(A_{n,J},\zeta)\gtrsim m(J_{k,n})$. Hence, using \eqref{eq:Jsum},
	\begin{equation}\label{eq:outer-tail}
		U_{q_n}(\zeta)
		\lesssim
		\sum_k
		\frac{1}{m(J_{k,n})^2}
		\sum_{J\in B(k,n)}m(J)^2
		\lesssim
		\varepsilon_n.
	\end{equation}
	Thus $\sum_{n\ge N}U_{q_n}(\zeta)<\infty$.

	It remains to consider the finitely many levels $n<N$. Fix such an $n$. If $\zeta\in K_n$, the preceding estimate applies. Otherwise, choose an open arc $V$ with $\zeta\in V$ and $\overline V\subset G_n$, and put $d=\dist(\overline V,K_n)>0$. If $J_{k,n}\cap V\neq\varnothing$, then \eqref{eq:whitney} implies $m(J_{k,n})\gtrsim d$. Since the Whitney arcs have pairwise disjoint interiors, only finitely many such $J_{k,n}$ occur. Each contains only finitely many arcs $J\in B(k,n)$, and every $A_{n,J}$ is a compact subset of $D$, whereas $\zeta\in E$. These terms contribute a finite amount. All remaining supports lie outside $V$, so their total contribution is bounded by $\dist(\zeta,\partial\D\setminus V)^{-2}\|q_n\|_1<\infty$. Thus $U_{q_n}(\zeta)<\infty$ for every $n<N$, and consequently $U_q(\zeta)<\infty$.

	Finally, define
	\begin{equation}\label{eq:outer-final-definition}
		O(z)
		=
		\exp\!\left(
		-\int_{\partial\D}
		\frac{\xi+z}{\xi-z}
		q(\xi)\,dm(\xi)
		\right).
	\end{equation}
	Then $O$ is outer and belongs to $\mathcal S$. By Lemma~\ref{lem:canonical-factor},
	$\mathcal F(O)=\{\zeta\in\partial\D:U_q(\zeta)<\infty\}=E$.
\end{proof}

\section{The dyadic construction}
For $k\ge1$, a binary word of length $k$ is an element of $\{0,1\}^k$, and will be denoted as
$w=\varepsilon_1\cdots\varepsilon_k\in\{0,1\}^k$. A binary word $w$ of length $k$ determines the integer $n(w)$ and the dyadic interval $J(w)$ defined by
\begin{equation}\label{eq:int1}
	n(w)=\sum_{j=1}^k \varepsilon_j 2^{k-j},
	\qquad
	J(w)=\left[\frac{n(w)}{2^k},\frac{n(w)+1}{2^k}\right)\subset[0,1).
\end{equation}
We use a word and the corresponding dyadic interval interchangeably when no confusion can arise.
If $\alpha=\sum_{i=1}^\infty \varepsilon_i2^{-i}$, where $\varepsilon_i=0$ or $1$, $i=1,2,\ldots$, write
$\alpha_k=\varepsilon_1\cdots\varepsilon_k$ for its prefix of length
$k$. For $n=1,2,\ldots$, let $\mathcal{D}_n=\{j 2^{-n}: j=0,1,2,\ldots, 2^n-1\}$, and $\mathcal D=\bigcup_n \mathcal{D}_n$ be the set of dyadic numbers in the interval $[0,1)$. Given $x\in\mathcal{D}$, the dyadic generation of $x$ is the minimal positive integer $n$ such that $x\in\mathcal{D}_n$. For dyadic numbers,
which have two eventually constant binary expansions, we use the expansion
ending in zeros. The main combinatorial statement is the following.
\begin{lemma}[Main lemma]\label{prop:familia}
	There exists a family $\mathcal{H}$ of words of finite length such that
	\begin{itemize}
		\item[(a)] A number $\alpha\in[0,1)$ is non-dyadic if and only if the set
		      \begingroup
		      \renewcommand{\theHequation}{H1}
		      \begin{equation}
			      H(\alpha)=\{k\geq1:\alpha_k\in\mathcal{H}\} \tag{H1} \label{H1}
		      \end{equation}
		      \endgroup
		      is infinite.
		\item[(b)] If $q$ is a dyadic number of generation $n$, then there exist exactly two dyadic intervals in $\mathcal{H}$ of generation greater than $n$ that have $q$ as an endpoint.
		\item[(c)] There exists an absolute constant $c_0>0$ such that if $q$ is a dyadic number of generation $n$, and $J=J(w)$, $w\in\mathcal{H}$, is a dyadic interval of generation $k>n$, not adjacent to $q$, then
		      \begingroup
		      \renewcommand{\theHequation}{H2}
		      \begin{equation}
			      2^{k}\lvert q-a(w)\rvert\geq c_0 2^{\frac{k-n}{2}}, \tag{H2} \label{H2}
		      \end{equation}
		      \endgroup
		      where $a(w)$ denotes the left endpoint of $J$.
	\end{itemize}
\end{lemma}
The rest of this section is devoted to constructing a family satisfying these requirements.
\begin{definition}[Maximal block decomposition]
	Every word $w\in\{0,1\}^k$ has a unique maximal block decomposition
	$w=b_1^{L_1}b_2^{L_2}\cdots b_m^{L_m}$, where
	$b_i\in\{0,1\}$, $b_i\ne b_{i+1}$, and each $b_i^{L_i}$ is a maximal
	block of consecutive equal digits. In particular, $L_1+\cdots+L_m=k$.
\end{definition}

\begin{definition}[Future after a block and heavy words]\label{def-heavy}
	Let $w=b_1^{L_1}\cdots b_m^{L_m}$ be its maximal block decomposition and
	assume $m>1$. For $i<m$, define the future of the $i$th block by
	$F_i(w)=L_{i+1}+\cdots+L_m$. The word $w$ is \emph{heavy} if
	$L_i\le F_i(w)$ for every $i<m$, with equality for at least one such $i$. All other words are \emph{light}; in particular, the constant words $0^k$ and $1^k$, $k=1,2,\ldots$, are declared light.
\end{definition}
For computations, the following equivalent formulation is useful.
\begin{lemma}[The operator $D$]\label{lem-D-equivalence}
	For a word $w=b_1^{L_1}\cdots b_m^{L_m}$, consider
	\begin{equation}
		D(w)=
		\begin{cases}
			\displaystyle \max_{1\leq i<m}\bigl((L_1+\cdots+L_i)+L_i\bigr), & m\geq2, \\[4pt]
			0,                                                              & m=1.
		\end{cases}
	\end{equation}
	Then $w$ is heavy if and only if $D(w)=|w|$, where $|w|$ denotes the length of $w$.
\end{lemma}
\begin{proof}
	We can assume $m>1$. Since $|w|=L_1+\cdots+L_m$, the condition $L_i\leq F_i(w)$ is equivalent to
	$(L_1+\cdots+L_i)+L_i\leq |w|$.
	Moreover, equality $L_i=F_i(w)$ holds for some $1\leq i<m$ if and only if $D(w)=|w|$.
\end{proof}
Let $\mathcal H_k$ denote the heavy words of length $k$. Set $\mathcal H=\bigcup_{k\ge1}\mathcal H_k$. The rest of this section is devoted to proving that this family
satisfies Lemma~\ref{prop:familia}.
\subsection{Condition (H1)}
\begin{proposition}\label{prop-finite-heavy}
	If $\alpha\in[0,1)$ is a dyadic point, then the set $H(\alpha)$ is finite.
\end{proposition}
\begin{proof}
	If $\alpha=0$ or $\alpha=1$, every prefix of $\alpha$ is light and so $H(\alpha)$ is empty. Assume now that $\alpha$ is dyadic. Then its binary expansion is eventually constant, that is, $\alpha$ can be written as
	\begin{equation}
		\alpha = b_1^{L_1} b_2^{L_2}\cdots b_{m-1}^{L_{m-1}} b_m^\infty,
		\label{eq:eventually-constant-run-decomp}
	\end{equation}
	where $b_j\in\{0,1\}$, $b_j\neq b_{j+1}$, $j=1,\ldots, m-1$ and $b_m^\infty$ denotes the infinite constant tail. We assume it is a block of zeros, but it does not matter in the argument. Set $P_j=L_1+\cdots+L_j$.
	For every $k>P_{m-1}$, the prefix $\alpha_k$ is
	\begin{equation}
		\alpha_k
		=
		b_1^{L_1} b_2^{L_2}\cdots b_{m-1}^{L_{m-1}}
		b_m^{\,k-P_{m-1}}.
		\label{eq:large-prefix-run-decomp}
	\end{equation}
	For every $k>P_{m-1}$,
	\begin{equation}
		D(\alpha_k)=\max_{1\leq i<m}(P_i+L_i).
		\label{eq:D-large-prefix-fixed-max}
	\end{equation}
	The right-hand side does not depend on $k$. Now
	$\alpha_k$ is heavy if and only if $D(\alpha_k)=k$.
	Hence among all prefixes with $k>P_{m-1}$, at most one value of $k$ satisfies that $\alpha_k$ is heavy. On the other hand, the remaining prefixes satisfy $k\leq P_{m-1}$. Therefore there are only finitely many heavy prefixes. Hence the set $H(\alpha)$ is finite.
\end{proof}
We next record the properties of $D$ needed below.
\begin{lemma}\label{lemma:ap1}
	\leavevmode\par
	\begin{enumerate}
		\item Let $w$ be a binary word that ends in the digit $b\in\{0,1\}$. Then,
		      \begin{equation}
			      D(w)=D(wb^K),
		      \end{equation}
		      for any $K\geq 1$.
		\item Let $w$ be a binary word of length $k$, final block of length $L_m$, and last digit $b\in\{0,1\}$. Appending $(1-b)$ creates a new block and
		      \begin{equation}
			      D(w(1-b))=\max\{D(w),k+L_m\}.
		      \end{equation}
	\end{enumerate}
	Moreover, if $\alpha\in[0,1)$ has binary expansion
	\begin{equation}
		\alpha=b_1^{L_1}b_2^{L_2}b_3^{L_3}\cdots,
		\qquad b_i\neq b_{i+1}, \ \ b_i\in\{0,1\}.
	\end{equation}
	If $L_{1}+\cdots+L_{r-1}<k\le L_1+\cdots+L_r$ and
	$R_r:=D(\alpha_k)$, then
	\begin{equation}
		R_r=\max\{R_{r-1},(L_1+\ldots+L_r)+L_r\}.
	\end{equation}
	Recall that $\alpha_k$ is the prefix of length $k$ of the number $\alpha$.
\end{lemma}
\begin{proof}
	$(1)$ The word $w$ ends with $b$, so adding $b^K$ does not produce more final blocks, so that $D$ remains the same.

	$(2)$ Note that adding $1-b$ means that the earlier final block is now a completed block. The old final block has length $L_m$. That is, we can now consider $d_m=k+L_m$ in the computation of $D(w(1-b))$. If $d_m< D(w)$, then $D(w(1-b))=D(w)$. Otherwise, $D(w(1-b))=k+L_m$. The formula follows from this observation.

	The recurrence follows by applying the second identity when a new block
	begins and the first identity while that block is extended.
\end{proof}
\begin{proposition}\label{prop-infinite-heavy}
	If $\alpha\in[0,1)$ is not a dyadic point, then $H(\alpha)$ is
	infinite.
\end{proposition}
\begin{proof}
	As $\alpha$ is not dyadic, its binary expansion is not eventually constant. Write the infinite block decomposition of $\alpha$ as
	\begin{equation}
		\alpha=b_1^{L_1}b_2^{L_2}b_3^{L_3}\cdots,
		\qquad b_i\neq b_{i+1}.
	\end{equation}
	Let $P_r=L_1+\cdots+L_r$. Consider integers $k$ with $P_{r-1}<k\leq P_r$. Hence the prefix $\alpha_k$ of $\alpha$ is heavy precisely when $P_i+L_i\leq k\quad\text{for every }i<r$,
	and equality holds for at least one $i<r$. Thus the only possible $k\in(P_{r-1},P_r]$ is
	\begin{equation}
		k=R_{r-1}:=\max_{1\leq i<r}(P_i+L_i).
	\end{equation}
	Set $S_r:=R_{r-1}-P_{r-1}$. Since
	$R_{r-1}\geq P_{r-1}+L_{r-1}$ and $L_{r-1}>0$, we have $S_r>0$.
	There exists $k\in(P_{r-1},P_r]$ such that $\alpha_k$ is heavy exactly
	when $S_r\leq L_r$. Moreover,
	\begin{equation}
		S_r-L_r=(R_{r-1}-P_{r-1})-L_r=R_{r-1}-(P_{r-1}+L_r)=R_{r-1}-P_r.
	\end{equation}

	We next derive a recurrence for $S_r$. An application of Lemma~\ref{lemma:ap1} yields
	\begin{equation}
		R_r=\max\{R_{r-1},P_r+L_r\}.
	\end{equation}
	Then,
	\begin{align}
		S_{r+1}
		 & =R_r-P_r
		=\max\{R_{r-1},P_r+L_r\}-P_r \\
		 & =\max\{R_{r-1}-P_r,L_r\}
		=\max\{S_r-L_r,L_r\}.
	\end{align}
	We argue by contradiction. Let us assume that $\alpha$ has only finitely many heavy prefixes. In that case, there exists $r_0\in\N$ such that every block $r\geq r_0$ contains no heavy prefix. Hence,
	\begin{equation}
		S_r> L_r, \ \ \ r\geq r_0.
	\end{equation}
	As $L_r>0$, we have $S_{r+1}=\max\{S_r-L_r,L_r\}< S_r, \ \ r\geq r_0$. Hence the positive integers $S_r$ form a strictly decreasing sequence for $r\geq r_0$, which is impossible, because such a sequence cannot remain positive indefinitely. Thus $H(\alpha)$ must contain infinitely many heavy prefixes.
\end{proof}
Propositions~\ref{prop-finite-heavy} and~\ref{prop-infinite-heavy} prove condition~\eqref{H1}.
\subsection{Adjacent intervals}
For a dyadic point, we first estimate the lengths of its heavy prefixes.
\begin{lemma}\label{lem-length-heavy}
	Let $\alpha$ be a dyadic number of generation $n$, with binary expansion $\alpha=\sum_{i=1}^{n}u_i 2^{-i}$, $u_i\in\{0,1\}$, $u_n=1$. Then $n<\max H(\alpha)\le 2n$. Furthermore, the prefix that attains the maximum is the only heavy prefix of length greater than $n$.
\end{lemma}
\begin{proof}
	Represent $\alpha$ as a binary word $u$, and write it in maximal block form
	\begin{equation}
		u=b_1^{L_1}b_2^{L_2}\cdots b_m^{L_m}.
	\end{equation}
	Set $E_i=L_1+\cdots+L_i$,  $1\le i\le m$. Since $u_n=1$, in the infinite expansion of $\alpha$, the finite word $u$ is followed by an infinite run of zeros. For every $k>n$, the prefix $\alpha_k$ has the form $\alpha_k=(u_1,u_2,\cdots ,u_n,0^{k-n})$. Thus all runs of $u$ are completed runs of $\alpha_k$, while the final
	zero-run $0^{k-n}$ is not counted. Hence, for every $k>n$,
	\begin{equation}
		D(\alpha_k)
		=\max_{1\le i\le m}(E_i+L_i).
	\end{equation}
	Define $R:=\max_{1\le i\le m}(E_i+L_i)$. Then $D(\alpha_k)=R$ for every $k>n$.
	Since $E_m=n$, we have $E_m+L_m=n+L_m$.
	Because $L_m\ge 1$, it follows that
	\begin{equation}
		R\ge E_m+L_m=n+L_m>n.
	\end{equation}
	Since $R>n$, we may take $k=R$ in the formula above. Therefore $D(\alpha_R)=R$.
	Hence $R\in H(\alpha)$. On the other hand, if $k>n$, then $D(\alpha_k)=R$. Thus $\alpha_k$ can be heavy only when $k=R$.
	All other possible heavy prefixes have length $k\le n$. Therefore
	\begin{equation}
		\max H(\alpha)=R.
	\end{equation}
	It remains to estimate $R$. For every $1\le i\le m$, we have $E_i\le n$ and $L_i\le n$. Therefore $E_i+L_i\le 2n$. Taking the maximum over $i$, we get $R\le 2n$. Combining this with $R>n$ and $\max H(\alpha)=R$, we obtain
	\begin{equation}
		n<\max H(\alpha)\le 2n,
	\end{equation}
	as we wanted to show. The proof shows that there can be only one heavy word
	of length greater than $n$, and it corresponds to the word that attains the
	maximum.
\end{proof}
\begin{lemma}\label{lem-adjacent}
	Let $q$ be a dyadic point of generation $n$. Let $A=A(q)$ be the set of dyadic intervals that have $q$ as an endpoint and are of generation greater than $n$. Then $A$ contains exactly two heavy intervals.
\end{lemma}
\begin{proof}
	This follows from the previous lemma.
\end{proof}
Lemma~\ref{lem-adjacent} proves part~(b) of Lemma~\ref{prop:familia}.
\subsection{Condition (H2)}
Condition (c) of Lemma~\ref{prop:familia} follows from the following result.
\begin{lemma}\label{lem-heavy-far}
	Let $q=J/2^n$ be a dyadic point of generation $n$, and let $k>n$ be an integer. If $w\in\{0,1\}^k$ is a heavy binary word and the dyadic interval it generates does not have $q$ as endpoint, then
	\begin{equation}
		2^k|q-\frac{n(w)}{2^k}|\geq \frac{1}{2}\cdot 2^{\tfrac{k-n}{2}}.
	\end{equation}
\end{lemma}
\begin{proof}
	Let $K=k-n$. Let $Q=J2^K$. If $2^k|q-\frac{n(w)}{2^k}|$ is already larger than
	$2^{K/2-1}$, there is nothing to prove. Thus we only need to rule out heavy words very
	close to $q$, but not adjacent to $q$.

	First suppose that the dyadic point obtained from the word $w$ lies to the right of $q$. Note that $q=J/2^n=J2^K/2^k=Q/2^k$, so $n(w)=Q+j$ with $j\geq1$; the case
	$j=0$ is the right-adjacent word. Let $\ell$ be the binary length of $j$, that is, the integer so that
	$2^{\ell-1}\leq j<2^\ell$.

	The last $K$ digits of $w$ are the $K$-digit binary expansion of $j$, hence begin with $K-\ell$ zeros. If it starts with a $1$, then the relative distance is already large, so that $K-\ell>0$.

	Thus $w$ has a completed block of
	zeros of length at least $K-\ell$, ending at position at least $k-\ell$. For this
	completed block, the quantity appearing in $D(w)$ is at least
	\begin{equation}
		D(w)\geq (k-\ell)+(K-\ell)=k+K-2\ell.
	\end{equation}
	If $\ell<K/2$, then $D(w)>k=|w|$, so $w$ is not heavy. Hence a heavy word must
	satisfy $\ell\geq K/2$, and therefore $j\geq 2^{K/2-1}$. Since
	\begin{equation}
		2^k|q-\frac{n(w)}{2^k}|=2^k|\frac{Q-(Q+j)}{2^k}|=|j|,
	\end{equation}
	the desired bound follows on the right with constant $1/2$. The left-hand case is analogous. \end{proof}

\section{Dyadic estimates and finite angular derivatives}\label{sec:dyadic-estimates}
Recall that $\mathcal D=\{j/2^k:k\ge1,\ 0\le j<2^k\}$ denotes the set of dyadic numbers. We first prove Theorem~\ref{thm:intro-countable} when $E=\{e^{2\pi ix}:x\in\mathcal{D}\}$. Recall that a word $w$ of length $k$ determines the dyadic interval $J(w)$
from \eqref{eq:int1}. For $w\in\mathcal H_k$, choose the two third-points
$t_{w,s}=(n(w)+s/3)/2^k$, $s=1,2$ in $J(w)$, and put $\delta_w=m(J(w))^2=2^{-2k}$.
We consider the sequence
\begin{equation}\label{eq:Z-dyadic}
	\mathcal Z=\{0\}\cup\bigl\{z_{w,s}:=(1-\delta_w)e^{2\pi i t_{w,s}}:
	w\in\mathcal H_k,\ k\ge1,\ s=1,2\bigr\}.
\end{equation}
It is clear that $\mathcal Z$ is a Blaschke sequence. Indeed, at generation $k$ there are $2^k$ dyadic intervals, and we place at most two zeros above each one. Hence there are at most $2^{k+1}$ zeros at generation $k$, and therefore $\sum_{z\in\mathcal Z}(1-|z|)\le 1+2\sum_{k\ge1}2^{-2k}2^k<\infty$.

Let $B$ be the Blaschke product with zero set $\mathcal Z$. We first verify that $B$ is an almost isometry.

\begin{lemma}\label{lem-separation}
	Let $q=J/2^n$ be a dyadic number of generation $n$ and let
	$\tau=(m+s/3)/2^k$, $s\in\{1,2\}$, where $m$ and $k$ are positive integers. Then
	$|q-\tau|\ge (3\cdot2^{\max\{n,k\}})^{-1}$.
\end{lemma}

\begin{proof}
	If $k\le n$, then
	$q-\tau=[3J-2^{n-k}(3m+s)]/(3\cdot2^n)$. The numerator is a nonzero integer,
	since the first term is divisible by $3$ and the second is not. Hence
	$|q-\tau|\ge(3\cdot2^n)^{-1}$. If $k>n$, then
	$q-\tau=[3J2^{k-n}-(3m+s)]/(3\cdot2^k)$, and the same divisibility argument
	gives $|q-\tau|\ge(3\cdot2^k)^{-1}$.
\end{proof}

We will repeatedly use the elementary estimate for $x,y\in\mathbb R$ such that $|x-y|\leq 1/2$,
\begin{equation}\label{basic-distance-estimate}
	|e^{2\pi i x}-(1-\delta)e^{2\pi i y}|^2\simeq \delta^2+|x-y|^2,
	\qquad 0<\delta<1.
\end{equation}
Let $z=(1-\delta)e^{2\pi i y}$ for some $0<\delta< 1$. If $I$ is an arc centred at
$e^{2\pi i c_I}$ and $z_I=(1-m(I))e^{2\pi i c_I}$, the same estimate yields
\begin{equation}\label{basic-carleson-estimate}
	|1-(1-\delta)e^{-2\pi i y}z_I|^2
	\simeq (\delta+m(I))^2+|y-c_I|^2.
\end{equation}

\begin{lemma}\label{lem-coarse}
	Let $I\subset[0,1)$ be a dyadic interval of generation $n$ and let $q$ be one of
	its endpoints. Then
	\begin{equation}
		\sum_{k\le n}\sum_{w\in\mathcal H_k}\sum_{s=1}^2\frac{1-|z_{w,s}|^2}{|e^{2\pi iq}-z_{w,s}|^2}\lesssim T_{\mathcal Z}(I).
	\end{equation}
\end{lemma}

\begin{proof}
	For $k\le n$ and $w\in\mathcal H_k$, Lemma~\ref{lem-separation} gives
	$|q-t_{w,s}|\ge2^{-n}/3$. If $c_I$ is the centre of $I$, then
	$|t_{w,s}-c_I|\le |t_{w,s}-q|+|q-c_I|\le |t_{w,s}-q|+m(I)/2
		\lesssim |t_{w,s}-q|$. Using \eqref{basic-distance-estimate} and
	\eqref{basic-carleson-estimate}, we obtain
	$|1-\overline{z_{w,s}}z_I|^2\lesssim|e^{2\pi iq}-z_{w,s}|^2$.
	Multiplying by $1-|z_{w,s}|^2$ and summing over $k\le n$ proves the claim.
\end{proof}

Let $q\in[0,1)$ be a dyadic point of generation $n$. For $k>n$ and
$w\in\mathcal H_k$, define the normalized distance
$M_{w,s}(q)=2^k|q-t_{w,s}|$, $s=1,2$. Since
$M_{w,s}(q)\asymp 2^k|q-n(w)/2^k|$, condition \eqref{H2} applies.

\begin{lemma}\label{lem-fine-reduction}
	Let $q$ be a dyadic point of generation $n$ and let $w$ be a heavy word of
	length $k>n$. Then
	\begin{equation}
		\frac{1-|z_{w,s}|^2}{|e^{2\pi iq}-z_{w,s}|^2}\lesssim M_{w,s}(q)^{-2}, s=1,2.
	\end{equation}
	In particular, there is an absolute constant $C$ such that, for every $k>n$,
	\begin{equation}
		\sum_{|w|=k}\sum_{s=1}^2M_{w,s}(q)^{-2}\le C.
	\end{equation}
\end{lemma}

\begin{proof}
	Write $h_k=2^{-k}$. Since $\delta_w=h_k^2$, we have
	$1-|z_{w,s}|^2\simeq h_k^2$. Also,
	$|e^{2\pi iq}-z_{w,s}|^2\simeq\delta_w^2+|q-t_{w,s}|^2
		\ge h_k^2M_{w,s}(q)^2$, which gives the first bound. Since $q$ is also a point of the dyadic grid of generation $k$, every dyadic interval of that generation at distance $j$ from $q$ has third-points with
	$M_{w,s}(q)\simeq j+1$. Hence
	$\sum_{|w|=k}\sum_{s=1}^2M_{w,s}(q)^{-2}
		\lesssim\sum_{j\ge0}(j+1)^{-2}<\infty$.
\end{proof}

\begin{lemma}\label{lem-heavy}
	There exists an absolute constant $C>0$ such that, for
	every dyadic point $q$ of generation $n$,
	\begin{equation}
		\sum_{k>n}\sum_{w\in\mathcal H_k}\sum_{s=1}^2M_{w,s}(q)^{-2}
		\le C.
	\end{equation}
\end{lemma}

\begin{proof}
	The contribution coming from adjacent chains is bounded because, by
	Lemma~\ref{lem-adjacent}, there are at most two heavy adjacent words, and
	for them $M_{w,s}(q)\in\{1/3,2/3\}$. For the non-adjacent contribution set
	$K=k-n$. By Lemma~\ref{lem-heavy-far}, $M_{w,s}(q)\gtrsim2^{K/2}$.
	Arguing as in Lemma~\ref{lem-fine-reduction}, we obtain
	$\sum_{w\in\mathcal H_k,\,w\notin A(q)}\sum_{s=1}^2M_{w,s}(q)^{-2}
		\lesssim\sum_{j\gtrsim2^{K/2}}j^{-2}\lesssim2^{-K/2}$. Summing over
	$k>n$ gives the claim.
\end{proof}

\begin{corollary}\label{cor-fine}
	Let $I$ be a dyadic interval of generation $n$, let $q$ be an endpoint of
	$I$, and put $\eta=e^{2\pi iq}$. Then
	\begin{equation}
		\sum_{k>n}\sum_{w\in\mathcal H_k}\sum_{s=1}^2
		\frac{1-|z_{w,s}|^2}{|\eta-z_{w,s}|^2}\le C,
	\end{equation}
	where $C$ is an absolute constant.
\end{corollary}

\begin{proof}
	Combine Lemmas~\ref{lem-fine-reduction} and~\ref{lem-heavy}.
\end{proof}

\begin{proposition}[Condition (AI1)]\label{prop-AI1}
	There exists an absolute constant $A_1>0$ such that every arc
	$I\subset\partial\D$ with $0<m(I)<1/2$ satisfies
	$\inf_{\xi\in I}S_{\mathcal Z}(\xi)\le A_1T_{\mathcal Z}(I)$.
\end{proposition}

\begin{proof}
	Let $I$ be a dyadic interval, let $q$ be one of its endpoints, and put
	$\eta=e^{2\pi iq}$. Lemma~\ref{lem-coarse} and Corollary~\ref{cor-fine}
	give $S_{\mathcal Z}(\eta)\le C_1T_{\mathcal Z}(I)+C_2$. The zero at the origin contributes exactly one to $T_{\mathcal Z}(I)$, and hence
	$S_{\mathcal Z}(\eta)\le(C_1+C_2)T_{\mathcal Z}(I)$. Thus \eqref{AI1}
	holds for dyadic arcs. For a general arc $I$, choose a dyadic subarc
	$J\subset I$ with $m(J)\geq m(I)/4$. Then
	$\inf_{\xi\in I}S_{\mathcal Z}(\xi)\le
		\inf_{\xi\in J}S_{\mathcal Z}(\xi)\lesssim T_{\mathcal Z}(J)$.
	Since $J\subset I$ and $m(J)\ge m(I)/4$, there is a constant $C>0$
	such that $|1-\overline a z_I|^2\le C^2|1-\overline a z_J|^2$ for every
	$a\in\D$. Therefore $T_{\mathcal Z}(J)\lesssim T_{\mathcal Z}(I)$, and
	\eqref{AI1} follows.
\end{proof}

\begin{lemma}\label{lem-kernel-grid}
	There exists an absolute constant $C>0$ such that, for every $\lambda>0$, $x\in[0,1]$, and $k\ge1$,
	\begin{equation}
		\sum_{|w|=k}\sum_{s=1}^2
		\frac{1}{\lambda^2+|x-t_{w,s}|^2}
		\leq C\left(\frac{1}{\lambda^2}+\frac{2^k}{\lambda}\right).
	\end{equation}
\end{lemma}

\begin{proof}
	The points $t_{w,s}$ are separated by at least $2^{-k}/3$. Grouping them according to their distance from $x$ gives
	$\sum_{|w|=k}\sum_{s=1}^2(\lambda^2+|x-t_{w,s}|^2)^{-1}
	\lesssim \lambda^{-2}+\sum_{j\ge1}(\lambda^2+j^22^{-2k})^{-1}$.
	Comparison with $2^k\int_0^\infty(\lambda^2+u^2)^{-1}\,du$ yields the result.
\end{proof}

\begin{proposition}[Condition (AI2)]\label{prop-AI2}
	There exists an absolute constant $A_2$ such that every arc $I\subset\partial\D$ with $0<m(I)<1/2$ satisfies $(1-|z_I|^2)T_{\mathcal Z}(I)\le A_2$.
\end{proposition}

\begin{proof}
	Write $\lambda=m(I)$ and $h_k=2^{-k}$. Then $1-|z_I|^2\simeq\lambda$, $1-|z_{w,s}|^2\lesssim h_k^2$, and by \eqref{basic-carleson-estimate},
	$|1-\overline{z_{w,s}}z_I|^2\gtrsim(\lambda+h_k^2)^2+|c_I-t_{w,s}|^2$. Lemma~\ref{lem-kernel-grid}, applied with $\lambda+h_k^2$, therefore gives
	\begin{equation}
		\lambda\sum_{w\in\mathcal H_k}\sum_{s=1}^2
		\frac{1-|z_{w,s}|^2}{|1-\overline{z_{w,s}}z_I|^2}
		\lesssim
		\min\left\{\frac{\lambda}{h_k^2},\frac{h_k^2}{\lambda}\right\}
		+\min\left\{\frac{\lambda}{h_k},h_k\right\}.
	\end{equation}
	Splitting at $h_k^2\simeq\lambda$, both sums over $k$ are geometric and uniformly bounded. The zero at the origin contributes at most one, so \eqref{AI2} follows.
\end{proof}

\begin{corollary}[Almost isometry]\label{cor-almost-isometry}
	The Blaschke product $B$ with zero set $\mathcal Z$ is an almost isometry, and $\mathcal F(B)=\{e^{2\pi i x}: x\in\mathcal{D}\}$.
\end{corollary}

\begin{proof}
	The fact that $B$ is an almost isometry follows from Propositions~\ref{prop-AI1} and~\ref{prop-AI2} together with the Garnett--Papadimitrakis criterion in Theorem~\ref{thm:GP-criterion}. We now show that $S_{\mathcal{Z}}(e^{2\pi i q})<+\infty$ if and only if $q\in\mathcal{D}$.

	Suppose first that $q$ is dyadic. Choose a dyadic interval $I$ having $q$
	as an endpoint. The proof of Proposition~\ref{prop-AI1} gives
	$S_{\mathcal Z}(e^{2\pi iq})\le C_1T_{\mathcal Z}(I)+C_2<\infty$.

	Now let $\alpha$ be non-dyadic. Its binary expansion is not eventually
	constant. By Proposition~\ref{prop-infinite-heavy}, there are infinitely
	many positive integers $k$ such that $w=\alpha_k$ is heavy. Let $J(w)$ be
	the dyadic interval generated by $w$. For such $k$,
	$\delta_w=2^{-2k}$. Moreover $\alpha\in J(w)$ and the two radial
	projections $t_{w,1},t_{w,2}$ also lie in $J(w)$, so
	$|\alpha-t_{w,s}|\le2^{-k}$, $s=1,2$. Using
	\eqref{basic-distance-estimate}, we obtain
	$|e^{2\pi i\alpha}-z_{w,s}|^2\lesssim2^{-2k}$, while
	$1-|z_{w,s}|^2\simeq\delta_w=2^{-2k}$. Thus each heavy prefix contributes
	a fixed positive amount:
	$(1-|z_{w,s}|^2)/|e^{2\pi i\alpha}-z_{w,s}|^2\gtrsim1$, $s=1,2$.
	Since there are infinitely many heavy prefixes, the full sum diverges.
	Characterization \eqref{eq:calF-Blaschke} now gives
	$\mathcal F(B)=\{e^{2\pi i x}:x\in\mathcal D\}$.
\end{proof}
We now transport the previous dyadic construction by a bi-Lipschitz homeomorphism of the circle. The next lemma is the matching result needed for this purpose.

\begin{lemma}[Bi-Lipschitz matching]\label{lem:matching}
	Let $A,B\subset\partial\D$ be countable dense sets. Then there exists a bi-Lipschitz
	homeomorphism $h:\partial\D\to\partial\D$ such that $h(A)=B$.
\end{lemma}

\begin{proof}
	Choose $a_*\in A$ and $b_*\in B$, cut the circles there, and write the corresponding parameter sets as the countable dense sets $D=\{t\in[0,1]:a_\ast e^{2\pi it}\in A\}$, $S=\{s\in[0,1]:b_\ast e^{2\pi i s}\in B\}$. Note that $0,1\in D\cap S$. We will construct an increasing bi-Lipschitz
	homeomorphism $f:[0,1]\to[0,1]$ with $f(D)=S$.

	Enumerate the interior points as $D\cap(0,1)=\{d_1,d_2,\ldots\}$ and
	$S\cap(0,1)=\{s_1,s_2,\ldots\}$. Pick two sequences $\{\varepsilon_n\},\{\eta_n\}$ of positive numbers such that $\sum_n\varepsilon_n<\infty$ and
	$\sum_n\eta_n<\infty$.

	Let $D_0=S_0=\{0,1\}$ and $f_0(x)=x$, $x\in[0,1]$. By induction, define an increasing homeomorphism $f_{n-1}:[0,1]\to[0,1]$ and $S_{n-1}\subset S$, $D_{n-1}\subset D$. Let $p_n$ be the first point in the
	enumeration of $D\setminus D_{n-1}$, and let $a_n<b_n$ be consecutive points of $D_{n-1}$
	with $p_n\in(a_n,b_n)$. Since $S$ is dense, choose
	$q_n\in (S\setminus S_{n-1})\cap(f_{n-1}(a_n),f_{n-1}(b_n))$ sufficiently close to
	$f_{n-1}(p_n)$. Define $g_n(p_n)= q_n$ and interpolate affinely on the two new subintervals. By taking $q_n$ close enough to $f_{n-1}(p_n)$, continuity ensures
	$\|g_n-f_{n-1}\|_\infty\le\eta_n/2$ and makes each new slope differ from
	the old one by a factor between $1-\varepsilon_n$ and
	$1+\varepsilon_n$.

	Next we use a backward step. Let $v_n$ be the first point of $S\setminus(S_{n-1}\cup\{q_n\})$, and put
	$x_n=g_n^{-1}(v_n)$. Choose a point $u_n\in D\setminus (D_{n-1}\cup\{p_n\})$ sufficiently close to
	$x_n$, inside the same component determined by the already matched points. Define $f_n(u_n)=v_n$ and interpolate affinely on the two new subintervals again. The resulting map $f_n$ satisfies $\|f_n-g_n\|_\infty\le\eta_n/2$, and its new slopes differ from
	those of $g_n$ by factors in the same interval. All previous
	correspondences remain fixed.

	Because the slope control is used twice at stage $n$, every slope of $f_n$
	lies between $m_n:=\prod_{j\le n}(1-\varepsilon_j)^2$ and
	$M_n:=\prod_{j\le n}(1+\varepsilon_j)^2$. Thus, with
	$m=\inf_nm_n>0$ and $M=\sup_nM_n<\infty$, we have
	$m|x-y|\le|f_n(x)-f_n(y)|\le M|x-y|$ uniformly in $n$. Moreover
	$\|f_n-f_{n-1}\|_\infty\le\eta_n$, so $(f_n)$ converges uniformly to an
	increasing bi-Lipschitz map $f:[0,1]\to[0,1]$, with $f(D)=S$.

	Finally define $h(a_*e^{2\pi it})=b_*e^{2\pi if(t)}$, $0\le t\le1$. It is a bi-Lipschitz homeomorphism of $\partial\D$ and $h(A)=B$.
\end{proof}

\begin{lemma}[Radial extension]\label{lem:radial}
	If $h:\partial\D\to\partial\D$ is bi-Lipschitz with constant $L\geq 1$, that is,
	\begin{equation}
		L^{-1}|\xi-\eta|\leq |h(\xi)-h(\eta)|\leq L|\xi-\eta|, \ \ \ \xi,\eta\in\partial\D,
	\end{equation} then the map defined by
	$H(0)=0$ and $H(r\xi)=rh(\xi)$ is a bi-Lipschitz homeomorphism of
	$\overline\D$ with the same constant.
\end{lemma}

\begin{proof}
	For $z=r\xi$ and $w=s\eta$, we have
	$|r\xi-s\eta|^2=(r-s)^2+rs|\xi-\eta|^2$ and
	$|H(z)-H(w)|^2=(r-s)^2+rs|h(\xi)-h(\eta)|^2$. Insert the two
	bi-Lipschitz inequalities for $h$ and use $L\ge1$.
\end{proof}

\begin{lemma}\label{lem:transport}
	Let $Z$ be a Blaschke sequence satisfying \eqref{AI1} and \eqref{AI2}. Let $h$ be a bi-Lipschitz
	homeomorphism of $\partial\D$, let $H$ be its radial extension, and
	$W=H(Z)$. Then $\mathcal F(B_W)=h(\mathcal F(B_Z))$, and $W$ satisfies
	\eqref{AI1}--\eqref{AI2}.
\end{lemma}

\begin{proof}
	Since $|H(z)|=|z|$, $z\in\D$, and $H$ is bi-Lipschitz,
	\begin{equation}
		\frac{1-|H(z)|^2}{|h(\xi)-H(z)|^2}\asymp
		\frac{1-|z|^2}{|\xi-z|^2}, \ \ z\in\D, \ \ \xi\in\partial\D\end{equation}
	so the Ahern--Clark criterion (\ref{eq:canonical-angular-criterion}) shows that $\mathcal F(B_W)=h(\mathcal F(B_Z))$.

	We now deduce that $W$ satisfies \eqref{AI1}--\eqref{AI2}. Let $I$ be an arc and write $J=h(I)$. Since $h$ is
	bi-Lipschitz, $m(J)\asymp m(I)$. For $z=r\zeta$, the elementary
	estimate $|1-\overline z z_I|\asymp(1-|z|)+m(I)+|\zeta-c_I|$ holds. Since $h(c_I)\in J$,
	we have $|h(c_I)-c_J|\lesssim m(J)$. Applying the bi-Lipschitz bounds to
	$h$ and $h^{-1}$, we obtain
	$m(J)+|h(\zeta)-c_J|\asymp m(I)+|\zeta-c_I|$. Since
	$H(z)=rh(\zeta)$, it follows that
	\begin{equation}\label{eq:Ttransport}
		T_W(J)\asymp T_Z(I).
	\end{equation}
	Since $S_W(h(\xi))\asymp S_Z(\xi), \xi\in I$, we deduce
	$\inf_{\eta\in J}S_W(\eta)\lesssim\inf_{\xi\in I}S_Z(\xi)
		\lesssim T_Z(I)\asymp T_W(J)$, which is \eqref{AI1}. Finally,
	$1-|z_I|^2\asymp m(I)\asymp m(J)\asymp1-|z_J|^2$. Together with
	\eqref{eq:Ttransport}, this transports \eqref{AI2} as well.
\end{proof}

We can now prove Theorem~\ref{thm:intro-countable}.
\begin{proof}[Proof of Theorem~\ref{thm:intro-countable}]
	Apply Lemma~\ref{lem:matching} to the set of dyadic points $\mathcal D$ and to the set $E$, and
	let $h$ be the corresponding bi-Lipschitz homeomorphism of $\partial\D$ with $h(\mathcal D)=E$. Let $H$ be the radial extension of $h$ and transport the
	zero set $Z$ of $B_{\mathcal D}$ by $H$.
	It is clear that $W=H(Z)$ is a Blaschke sequence. Lemma~\ref{lem:transport} gives that the Blaschke product $B_W$ with zeros $W$ is an almost isometry and
	$\mathcal F(B_W)=h(\mathcal F(B_{\mathcal D}))=h(\mathcal D)=E$.
\end{proof}

\printbibliography[title={References}]

@book{Shapiro,
  author    = {Shapiro, Joel H.},
  title     = {Composition Operators and Classical Function Theory},
  publisher = {Springer-Verlag},
  address   = {New York},
  year      = {1993}
}

@article{AhernClark,
  author  = {Ahern, Patrick R. and Clark, Donald N.},
  title   = {Radial limits and invariant subspaces},
  journal = {American Journal of Mathematics},
  volume  = {92},
  year    = {1970},
  pages   = {332--342}
}

@book{Mashreghi,
  author    = {Mashreghi, Javad},
  title     = {Derivatives of Inner Functions},
  series    = {Fields Institute Monographs},
  volume    = {31},
  publisher = {Springer},
  address   = {New York},
  year      = {2013}
}

@book{Duren,
  author    = {Duren, Peter L.},
  title     = {Theory of {$H^p$} Spaces},
  publisher = {Academic Press},
  address   = {New York},
  year      = {1970}
}

@book{CMR,
  author    = {Cima, Joseph A. and Matheson, Alec L. and Ross, William T.},
  title     = {The Cauchy Transform},
  series    = {Mathematical Surveys and Monographs},
  volume    = {125},
  publisher = {American Mathematical Society},
  address   = {Providence, RI},
  year      = {2006}
}

@book{Garnett,
  author    = {Garnett, John B.},
  title     = {Bounded Analytic Functions},
  edition   = {Revised First Edition},
  publisher = {Springer},
  address   = {New York},
  year      = {2007}
}

@article{GP91,
  author  = {Garnett, John B. and Papadimitrakis, Michael},
  title   = {Almost isometric maps of the hyperbolic plane},
  journal = {Journal of the London Mathematical Society. Second Series},
  volume  = {43},
  year    = {1991},
  pages   = {269--282}
}

@article{IvriiNicolau2024,
  author  = {Ivrii, Oleg and Nicolau, Artur},
  title   = {Analytic mappings of the unit disk which almost preserve hyperbolic area},
  journal = {Proceedings of the London Mathematical Society},
  volume  = {129},
  number  = {5},
  year    = {2024},
  pages   = {e70001},
  doi     = {10.1112/plms.70001}
}

@misc{IvriiNicolau2025,
  author      = {Ivrii, Oleg and Nicolau, Artur},
  title       = {Analytic mappings of the unit disk with bounded compression},
  year        = {2025},
  eprint      = {2507.15200},
  eprinttype  = {arxiv},
  eprintclass = {math.CV}
}

\vspace{-0.25\baselineskip}
{\footnotesize\noindent Departament de Matemàtiques, Universitat Autònoma de Barcelona, 08193 Bellaterra, Barcelona, Spain. \textit{Contact:} \href{mailto:angel.lorenzo@uab.cat}{\texttt{angel.lorenzo@uab.cat}}\par}

\end{document}